\documentclass{article}

\usepackage[margin=1.5in]{geometry} 
\usepackage[utf8]{inputenc}
\usepackage{amsmath}
\usepackage{amsthm}
\usepackage{amsfonts}
\usepackage{amssymb}
\usepackage{graphicx}
\usepackage[]{algorithm2e}
\usepackage{subfig}
\usepackage[utf8]{inputenc}
\title{A Peculiarity in the Parity of Primes}
\author{Debayan Gupta
\\MIT
\\debayan@mit.edu
\and
Mayuri Sridhar
\\MIT
\\mayuri@mit.edu
}
\date{\today}

\begin{document}

\maketitle


\begin{abstract}

\noindent

We create a simple test for distinguishing between sets of primes and random numbers using just the {\sc sum-of-digits} function.

We find that the sum-of-the-digits of prime numbers does not have an equal probability of being odd or even. The authors know of no reason why prime numbers should bias themselves towards a particular parity in their sums of digits, but our empirical tests show a very strong bias; strong enough that we are able to devise a test to reliably differentiate between collections of prime numbers versus random numbers by looking {\it only} at their sums of digits. We are also able to create similar tests for {\it products} of primes. We are even able to test ``tainted'' sets with mixtures of primes and random numbers: as the percentage of (randomly chosen) prime numbers in a set of random numbers is varied, we get a reliable, linear change in our parity measure.

For example, when we add up the digits of prime numbers in base 10, their sum is significantly more likely to be odd than even. This effect persists across base changes, although {\it which} parity is more common might change. Note that the last digit being odd in base 10 simply reverses the parity. We have tested this for the first fifty million primes -- not primes up to 50,000,000, but the first 50,000,000 prime numbers -- and have found that this effect persists, and does so in a {\it predictable} manner. The effect is quite significant; for 50,000,000 primes in base 10, the number of primes which have an odd {\sc sum-of-digits} is about an order of magnitude farther away from the mean than expected.

We have run multiple tests to try and understand the source of this bias, including investigating primes modulo random numbers and adjusting for Chebyshev's bias. None of these tests yielded any satisfactory explanation for this phenomenon.

\end{abstract}
\newpage


\section{Introduction}

Conventional wisdom says that prime numbers should look superficially ``random-ish'' -- unless one actually checks for primality, there should not be any easy way to differentiate between a collection of random numbers and a collection of primes. There are some basic tests, like the fact that prime numbers (other than 2) are odd, or more complicated ones, such as the Lucas–Lehmer primality test for Mersenne numbers, or the AKS test~\cite{agrawal2004primes}. However, we understand why these tests work, i.e., which underlying properties of primes are exploited for these tests. 

In this paper, we explore a function which should not have anything to do with primality: the parity of the {\sc sum-of-digits}. The authors know of no reason why prime numbers should bias themselves towards a particular parity in their sums of digits, but our empirical tests show a very strong bias; strong enough that we are able to devise a test to reliably differentiate between collections of prime numbers versus random numbers by looking {\it only} at their sums of digits.


Our test involves {\it no} computation beyond addition, yet is capable of acting as a very strong judge, differentiating between groups of primes and random integers, even if the last digit of each number is removed. For a set of positive integers, $S$, we can test how likely they are to be a set of primes (versus random numbers)\footnote{We look at ``tainted'' sets of random numbers mixed with varying percentages of primes in section~\ref{sec:mixed}.} simply by summing up the digits of each element in $S$.  We denote the set of the ``sums of digits'' of the elements in $S$ as $D$. If the elements of $S$ were random integers, we would expect that approximately half the values in $D$ would be even and half would be odd. Conventional wisdom suggests that this should also be true if $S$ were a set of (randomly chosen) primes.

There are some papers which perform theoretical analyses on the subject~\cite{mauduit2010probleme,drmota2009primes}, and reach the conclusion that the
sums of digits of prime numbers should be balanced between odd and even paritites. When we compute the parity of the {\sc sum-of-digits} for a set of random numbers, we find that the results are within a few standard deviations from the mean. In fact, they are within one standard deviation almost all of the time. However, in the case of prime numbers, this does not hold true -- the results are {\bf significantly} farther from the mean than expected.

When we first looked at this, our suspicion was that we were just looking at some sort of statistical artifact -- a quirk that would be smoothed out once we considered larger samples. As we found out, this is not the case. In fact, {\bf the effect increases at a predictable rate}. In the first 50,000,000 prime numbers, we find 25,032,384 primes with an odd {\sc sum-of-digits}. This implies that there are 32384 more primes with an odd {\sc sum-of-digits} than we would expect. In our first round of analysis, we compared this with random numbers. We assumed (and confirmed, empirically) that random numbers have a sum of digits which has a 50\% chance of being odd or even. We expected this to behave like a binomial distribution with $p=0.5$ and $n=5\times10^7$. For this setting, we expected the set of numbers to have a mean of $np = 2.5\times10^7$ and a standard deviation of $\sqrt{np(1-p)} = 3535.53$.

This suggests that prime numbers do not behave like random numbers with respect to this metric. Thus, our assumption that ``primes superficially look like random numbers'' was incorrect.
In the next sections of this paper, we outline some more detailed tests we ran to compare prime numbers and random numbers with respect to our metric.

\begin{itemize}

\item Section~\ref{sec:parity} discusses how we compare sets of primes and random numbers, as well as expected values, our null hypothesis, and the likeliness of our empirical results.

\item Section~\ref{sec:testset} describes our algorithm for differentiating between sets of primes and random numbers, finds a simple quadratic curve that can predict the z-scores of sum-of-digits-parity for sets of primes of various sizes (which means we can predict the deviation from the mean for other sample sizes), and calculates our probability of success for various sample sizes.

\item Section~\ref{sec:product} looks at products of two primes and finds that we can still differentiate between these and products of random numbers, even if the random numbers are artificially biased to mimic primes.

\item Section~\ref{sec:mixed} lists our observations for sets of random number ``tainted'' with varying amounts of primes and observes a linear change in our metric as we vary the percentage of primes. This further solidifies our previous hypothesis that this is a predictable effect specific to primes.

\item Section~\ref{sec:additional} discusses a number of other tests we performed, such as adjusting for Chebyshev's bias, to see if our sum-of-digits effect was simply a symptom produced by one of these known biases. We did not find any satisfactory explanation.

\item Section~\ref{sec:conclusion} concludes and discusses future directions.

\end{itemize}


\section{Primes Parity Testing}\label{sec:parity}

In this section, we shall describe our process of comparing sets of primes and random numbers, and note the (significant) disparities we observe.

First, we choose sample sets of prime numbers (uniformly at random from the first 50 million primes), with sample sizes
varying from 100,000 to 5,000,000. We then test the hypothesis that these prime numbers behave the same way as random numbers: If our hypothesis is true,
then we would expect approximately half of the sample size to have an even {\sc sum-of-digits}. 

For a given trial, we draw a sample of size $s$ from the first 50,000,000 prime numbers. We can define a random variable $X_i$ as the number of samples
that have an even {\sc sum-of-digits} in trial $i$. If our hypothesis is true, then, $X_i$ is a binomial random variable, with expectation and variance as defined below:
\begin{align*}
\mathbb{E}[X_i] &= \frac{s}{2},\\
Var[X_i] &= \frac{s}{4}\,.
\end{align*}

To reduce the variance, for a fixed sample size, we ran $t$ trials. Thus, we can define the variable $\overline{X}$ as the average value of $X_i$ over all
$t$ trials. In particular, $\overline{X}$ has expectation and variance of
\begin{align*}
\mathbb{E}[\overline{X}] &= \frac{s}{2},\\
Var[\overline{X}] &= \frac{s}{4t}\,.
\end{align*}

In our experiments, we chose $t=1000$. As a baseline for our results, we also sampled random numbers. In particular, for a fixed sample size $s$, we sampled
uniformly at random from odd numbers in the range [3, m], where $m$ is the value of our largest prime number. This ensured that we were sampling numbers
from approximately the same range. Our results are summarized in Table \ref{table:full_results}.

\begin{table}[!htbp]
\centering
\begin{tabular}{| c | c | c | } 
 \hline
 Sample Size &
 \multicolumn{1}{|p{4cm}|}{\centering Average Number of Samples with Even Parities in Primes} &
 \multicolumn{1}{|p{4cm}|}{\centering Average Number of Samples with Even Parities in Random Numbers}\\
 \hline\hline
	100,000 & 49,936 & 49,995\\
	200,000 & 99,871 & 99,997\\
	300,000 & 149,803 & 149,998\\
	400,000 & 199,731 & 199,990\\
	500,000 & 249,666 & 249,989\\
	600,000 & 299,617 & 300,005\\
	700,000 & 349,566 & 349,983\\
	800,000 & 399,476 & 400,012\\
	900,000 & 449,415 & 450,006\\
	1,000,000 & 499,320 & 499,991\\
	2,000,000 & 998,707 & 1,000,024\\
	3,000,000 & 1.498,073 & 1,499,951\\
	4,000,000 & 1,997,360 & 1,999,966\\
	5,000,000 & 2,496,799 & 2,499,976\\[1ex] 
 \hline
\end{tabular}
\caption{$\overline{X}$ over 1,000 trials, for varying sample sizes when drawn from the prime numbers
distribution and when drawn from the random numbers distribution.}
\label{table:full_results}
\end{table}

We can normalize these results by subtracting the mean and dividing by the standard deviations.
We take the absolute value of this quantity to get the average number of standard deviations that the empirical data varies from the mean.
These results are summarized in Table \ref{table:normalized_results}.

\begin{table}[!htbp]
\centering
\begin{tabular}{| c | c | c | } 
 \hline
 Sample Size &
 \multicolumn{1}{|p{4cm}|}{\centering Z-Score of $\overline{X}$ for Primes} &
 \multicolumn{1}{|p{4cm}|}{\centering Z-Score of $\overline{X}$ for Random Numbers}\\
 \hline\hline
	100,000 & 12.80 & 1.10\\
	200,000 & 18.24 & 0.42\\
	300,000 & 22.79 & 0.20\\
	400,000 & 26.86 & 1.04\\
	500,000 & 29.85 & 0.96\\
	600,000 & 31.24 & 0.44\\
	700,000 & 32.80 & 1.31\\
	800,000 & 37.07 & 0.86\\
	900,000 & 38.99 & 0.43\\
	1,000,000 & 43.03 & 0.56\\
	2,000,000 & 57.83 & 1.07\\
	3,000,000 & 70.36 & 1.78\\
	4,000,000 & 83.49 & 1.07\\
	5,000,000 & 90.53 & 0.69\\[1ex] 
 \hline
\end{tabular}
\caption{The number of standard deviations from the mean for $\overline{X}$ when drawn from the prime numbers distribution and when drawn from the random numbers distribution,
over 1,000 trials.}
\label{table:normalized_results}
\end{table}

As we can see from the data above, for random numbers, the Z-score for $\overline{X}$ stays low, as expected. For prime numbers, however, even at a relatively small sample size,
like 100,000 numbers, there is a significant deviation from the expected value.

We can apply Chebyshev's bound using the Z-scores for the primes to calculate the probability that we will see this deviation, assuming that our null hypothesis holds.
These results are summarized in Table \ref{table:chebyshev}.

\begin{table}[!htbp]
\centering
\begin{tabular}{| c | c | c | } 
 \hline
 Sample Size &
 \multicolumn{1}{|p{4cm}|}{\centering Max Probability of Drawing Empirical Data}\\
 \hline\hline
	100,000 & $6.1 \times 10^{-3}$\\
	200,000 & $3.0 \times 10^{-3}$\\
	300,000 & $1.9 \times 10^{-3}$\\
	400,000 & $1.4 \times 10^{-3}$\\
	500,000 & $1.1 \times 10^{-3}$\\
	600,000 & $1.0 \times 10^{-3}$\\
	700,000 & $9.3 \times 10^{-4}$\\
	800,000 & $7.3 \times 10^{-4}$\\
	900,000 & $6.6 \times 10^{-4}$\\
	1,000,000 & $5.4 \times 10^{-4}$\\
	2,000,000 & $3.0 \times 10^{-4}$\\
	3,000,000 & $2.0 \times 10^{-4}$\\
	4,000,000 & $1.4 \times 10^{-4}$\\
	5,000,000 & $1.2 \times 10^{-4}$\\[1ex] 
 \hline
\end{tabular}
\caption{Chebyshev's bound on drawing the empirical values of $\overline{X}$ from the prime number distribution, for varying sample sizes over 1,000 trials, assuming that our null hypothesis holds.}
\label{table:chebyshev}
\end{table}

As shown from the above results, the probability of generating our empirical data assuming our null hypothesis holds is quite low.
Thus, we can choose a significance level of 0.01 and reject our null hypothesis that prime numbers and random numbers follow the same distribution,
when looking at the criteria of the parity of the {\sc sum-of-digits}.


\section{Testing for Sets of Prime Numbers}\label{sec:testset}

In the previous section, we showed that sets of prime numbers and random numbers behave differently when it comes to the {\sc sum-of-digits} metric. In this section, we shall exploit this difference in behavior to build a simple algorithm that we can use to distinguish between sets of prime numbers and sets of random numbers. We note that our algorithm is simpler than previous testing algorithms, like tests based on prime density. Moreover, it's easily generalizable, since it doesn't rely on the range of the prime numbers produced. For simplicity, we will focus on the case of distinguishing between a set of random odd numbers and a set of prime numbers.

\begin{figure}[!htbp]
\centering
 \includegraphics[width=10cm, height=6cm,keepaspectratio]{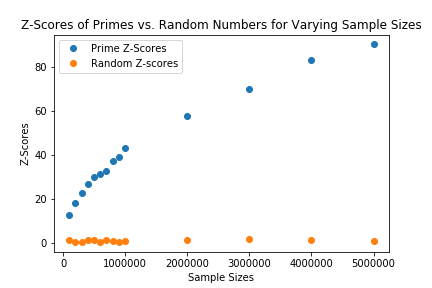}

\caption{Z-scores for primes and random numbers, with varying sample sizes.}
\label{fig:z_by_sample}
\end{figure}

To design our algorithm, we first plot the Z-score for $\overline{X}$, with respect to sample size for prime numbers and random numbers, as seen in Figure \ref{fig:z_by_sample}.

As we can see in Figure \ref{fig:z_by_sample}, the Z-score for random numbers
is approximately constant, with respect to sample size. However, the Z-score for prime numbers increases with our sample size dramatically.
In fact, we note that the Z-scores increase with the sample size at a predictable rate.

\begin{figure}[!htbp]
\centering
 \includegraphics[width=10cm, height=6cm,keepaspectratio]{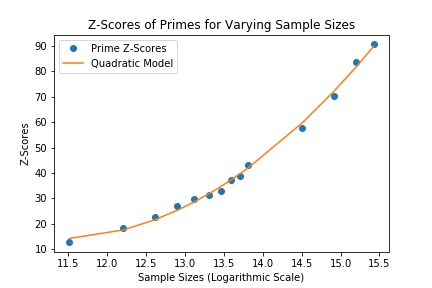}

\caption{Z-scores of primes with varying sample sizes and quadratic model fit. We note that the sample sizes here are on a logarithmic scale.}
\label{fig:quad_model}
\end{figure}

To analyze this, as seen in Figure \ref{fig:quad_model}, we tried to fit a quadratic model to our data. The model $$4.49 t^2 -101.62 t +  589.13$$
with $t = \ln x$ proved to be the model of best fit, with a sum-of-squares deviation of 23.33.
We note that curve fitting with $\ln(x)$ provided a significantly better fit than curve fitting with $x$. We suspect this is partially because the coefficients were so small (on the order of $10^{-12}$ for the second order term); this may have resulted in numerical issues.

The z-scores for random numbers are linear, with no dependence on the sample size, and never more than 2 standard deviations away from expectation. This suggests that Z-scores could be a useful distinguisher between sets of randomly drawn prime numbers and random numbers.

Formally, let's assume we have a source $G$ which can produce samples from some distribution $\mathbb{E}$.
We would like to determine whether $\mathbb{E}$ is the distribution of prime numbers or the distribution of random numbers. To do this,
we first define a simple helper function, as described in Algorithm \ref{alg:helper}.

\begin{center}
\begin{algorithm}[H]
 \KwData{A set $S$ of numbers, which are either prime or random}
 \KwResult{The number of cases with even parity}
 num\_even = 0\;
 \For{number in S}{
  Calculate the parity of the sum of digits $p$\;
  \If{$p$ is even}{
  	num\_even ++\;
   }
 }
 return num\_even

 \caption{Simple helper function which computes the number of numbers in $S$ with an even {\sc sum-of-digits}.}
 \label{alg:helper}
\end{algorithm}
\end{center}

This helper function, for a given set of integers $S$, simply computes the number of elements in $S$ with an even {\sc sum-of-digits}, which we can use in our testing
algorithm. Thus, we can design the following testing algorithm, described in Algorithm \ref{alg:tester}.

\begin{center}
\begin{algorithm}[H]
 \KwData{A source $G$ which produces samples from either the prime number distribution or a random number distribution.}
 \KwResult{``Yes'' if $G$ is sampling from the prime number distribution, with high probability.}
 sample\_size = $10^5$\;
 \For{trial in $1 \ldots 1000$}{
  Draw a sample $s$ of size sample\_size from $G$\;
  Run Algorithm \ref{alg:helper} on $s$ and let $p_t$ be the result\;
  }
 Compute the average value of $p_t$, $p_{avg}$ over the trials\;
  \BlankLine
 $exp_{avg} = \frac{\text{sample\_size}}{2.}$\;
  \BlankLine
 $std_{avg} = \sqrt{\frac{\text{sample\_size}}{4*1000}}$\;
  \BlankLine
 $z_{avg} = \frac{|p_{avg} - exp_{avg}|}{std_{avg}}$\;
   \If{$z_{avg} > 5$}{
  	Return ``Yes''\;
   }
   Return ``No''\;

 \caption{Simple testing algorithm, which differentiates between primes and random numbers, using the parity of the sums of digits.}
 \label{alg:tester}
\end{algorithm}
\end{center}

For simplicity, we chose fixed sample sizes of 100,000, 1,000 trials, and a threshold of 5 standard deviations. These can vary, based on $G$.

We can analyze the behavior of our algorithm for these settings. In particular, we assume that $G$ samples uniformly at random from either
the prime distribution or the random number distribution. We analyze the behavior of $G$, assuming that it samples uniformly at random from the first 50 million
primes, if sampling from the prime distribution. Based on our earlier results, we expect our algorithms behavior to improve when sampling from a larger range.

For any given draw from the prime number distribution, we know that we draw a prime with an even {\sc sum-of-digits} with probability $\frac{24967616}{50\times10^{6}}$, or $q=0.49935232$.
For $z_{avg}$ to be greater than 5, we require:
$$\frac{|p_{avg} - \frac{\text{sample\_size}}{2}|}{std_{avg}} > 5.$$

For simplicity, since $q < 0.5$, we only analyze the case where $z_{avg}$ is at least 5 standard deviations below the mean. This implies that for $z_{avg}$ to be greater than 5, the following must hold:

\begin{align*}
\frac{ \frac{\text{sample\_size}}{2} - p_{avg}}{std_{avg}} &> 5\,\\
\frac{\frac{\text{sample\_size}}{2} - p_{avg}}{\sqrt{\frac{\text{sample\_size}}{4000}}} &> 5\,\\
\frac{\text{sample\_size}}{2} - p_{avg} &> 5\sqrt{\frac{\text{sample\_size}}{4000}}\,\\
p_{avg} &< \frac{\text{sample\_size}}{2} - 5\sqrt{\frac{\text{sample\_size}}{4000}}.
\end{align*}

For $p_{avg}$ to satisfy our bound, we know that the sum of $p_t$ must satisfy
$$\sum_{t=1}^{1000} p_t < 1000(\frac{\text{sample\_size}}{2} - 5\sqrt{\frac{\text{sample\_size}}{4000}}).$$

However, the sum of $p_t$ is simply a binomial random variable (when $\mathbb{E}$ is the distribution of the first $5 \times 10^7$ primes),
with a success probability of $p=0.49935232$ and $n=1000*\text{sample\_size}$. Thus, we can calculate
the binomial cumulative distribution function.

With our settings of sample\_size = $10^5$ and $t=1000$, 5 standard deviations below the expected value becomes $(49975 \times 10^3)$. We want to compute the binomial CDF
of getting less than $(49975 \times 10^3)$ numbers with an even {\sc sum-of-digits}, in $10^8$ trials, where each trial has a success probability of $q$.
This happens with high probability - about $(1-1.89\times 10^{-16})$ of the time, the sum of $p_t$ will be lower than 5 standard deviations below the expected value. Thus, when the input is prime, our algorithm returns ``Yes'' with high probability. When the input is a set of random numbers, our algorithm
will return incorrectly if and only if they deviate from the mean by 5 standard deviations. Again, for our suggested values, this happens with probability $2.87\times 10^{-7}$. Thus, with high probability,
our algorithm can distinguish between prime numbers and random numbers, given enough\footnote{Just how many samples are ``enough'' is an area for future work. We have been able to reliably judge sample sizes in the thousands, but perhaps there are better scoring mechanisms than Z-scores.} samples.


\section{Prime Product Testing}\label{sec:product}

Given our observations with prime numbers, we decided to look at what happened with products of primes. That is, we wanted to see if there were any patterns in numbers of the form $pq$ where $p$ and $q$ were both prime numbers. 

From our previous trials, we know that a set of randomly chosen odd numbers behaves like a binomial random variable, with expected value $n/2$ and standard deviation $\sqrt{n/4}$. We also know that a set of randomly chosen prime numbers does not behave this way. 

We ran a similar testing framework to investigate how the product of two prime numbers behaved. If it behaves like a random number, then, for any given sample size, $s$, for any given trial,
we expect the parity to be distributed like a binomial random variable $X_i$, which has $\mathbb{E}[X_i] = \frac{s}{2}$ and $Var[X_i] = \frac{s}{4}$. Thus, we can look at the average $X_i$ value,
over 100 trials, which we denote as the random variable $\overline{X}$.

If the products behave like random numbers, then
$$\mathbb{E}[\overline{X}] = \frac{s}{2}$$
$$Var[\overline{X}] = \frac{s}{4 \times 100}$$

\begin{figure}[!htbp]
\centering
 \includegraphics[width=10cm, height=6cm,keepaspectratio]{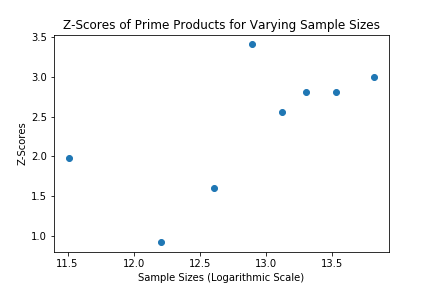}

\caption{Z-scores for prime products, with varying sample sizes. Note that
the sample size is plotted on a logarithmic scale.}
\label{fig:prod_by_sample}
\end{figure}

Thus, from here, we can calculate statistics about the prime products, using the empirical samples. Our normalized Z-scores for prime products, for varying sample sizes are plotted
in Figure \ref{fig:prod_by_sample}.

We can see that, like prime numbers, the distribution of prime products doesn't follow the same distribution as random numbers. The major observations we made were that the Z-scores
we computed for the products were often above 1 standard deviation. While this isn't as strong a result as we showed for prime numbers, this implies that we can design a simple
distinguishing algorithm which distinguishes between prime numbers and random numbers with a success probability significantly greater than 50\%.

Moreover, we note that the Z-scores for prime products are also dependent on sample size. This implies that if we test with larger sample sizes - for instance, if $G$ could produce
infinitely many samples - we could use this metric to distinguish between prime products and random numbers, as well.

We hypothesized that, {\it perhaps, the difference in parities in primes was causing the difference in parities in their products}. I.e., perhaps what we were seeing is what always happens when you look at a set of products of two numbers where the two numbers were drawn from an already-biased source. In that case, what we were observing for prime-products would be merely a symptom of the bias in the primes and not something deeper.

\begin{figure}[!htbp]
\centering
 \includegraphics[width=10cm, height=6cm,keepaspectratio]{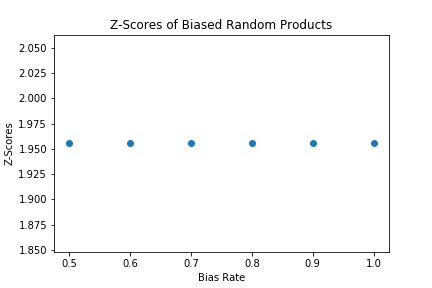}

\caption{Z-scores for biased random numbers, with varying bias rates and $s=400,000$.}
\label{fig:biased_rand}
\end{figure}

To test this, we artificially created sets of random numbers with varying
levels of parity bias. That is, we choose a sampling rate $r$, where $r$ ranged from 0.5 to 1.0, with a fixed sample size of $s=400,000$. For each trial, we generated $3s$ random numbers,
which we separated into two different lists $s_{even}$ and $s_{odd}$, where the former contained just elements with a sum of digits with even parity and the latter contained the rest.
We chose $rs$ numbers from $s_{even}$ and the remaining $(1-r)s$ numbers from $s_{odd}$. Thus, like prime numbers, our samples were
biased in the parity of their sum of digits. From here, we randomly sampled 2 numbers from our biased list, and calculated the product $s$ times, for a single trial. We ran 100 trials,
for each sampling rate $r$ and computed the average Z-scores. Our results are plotted in Figure \ref{fig:biased_rand}. We note that the Z-scores are represented in absolute-value form; we also note that, in all cases, our random products are slightly more likely to have an {\it even} sum of digits parity.

These results for random products vary significantly from the results we found for the prime products. In particular, we still see some variation from the expected value, with relatively high z-scores, however, the z-scores are constant, with respect to the sampling rate, rather than increasing. Notably, even the random products with $r=0.5$ show high z-scores. We hypothesize that this may be due to the fact that all the products are ``more composite'' by construction. Since the primes are biased towards having an odd sum of digits parity, the random products are biased towards having an even sum of digits parity. 

This suggests that {\it increasing the bias in the parity of the component random numbers does not increase the bias in the parity of the products}.
Thus, we can conclude that the bias in the primes is not the cause of the bias in the prime products; rather, {\bf the bias
in the prime products is due to some intrinsic property of prime numbers}.


\section{Prime and Random Mixed Testing}\label{sec:mixed}

The final set of tests that we explored for distinguishing between prime numbers and random numbers is to test whether we could use the Z-score as a metric
to approximate the percentage of primes that were in a set $S$.

\begin{figure}[!htbp]
\centering
 \includegraphics[width=10cm, height=6cm,keepaspectratio]{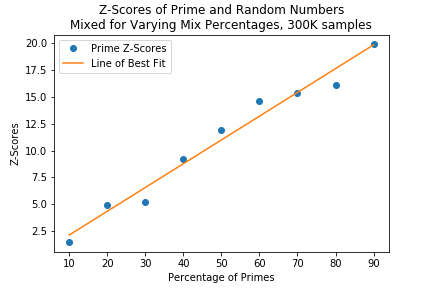}

\caption{Z-scores for prime and random numbers mixed. The percentage of prime numbers
is varied between 10\% and 90\%. The sample size is fixed at 300,000 samples.}
\label{fig:prime_ran;ss:300k}
\end{figure}

\begin{figure}[!htbp]
    \centering
    \subfloat[Sample Size of 100,000]{{\includegraphics[width=5.5cm]{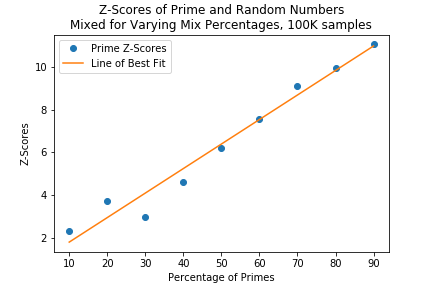} }}%
    \qquad
    \subfloat[Sample Size of 500,000]{{\includegraphics[width=5.5cm]{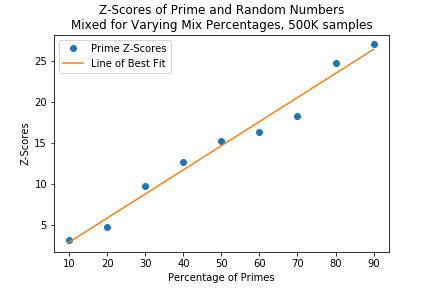} }}%
\caption{Z-scores for prime and random numbers mixed. The percentage of prime numbers
is varied between 10\% and 90\%. The sample size is fixed at 100,000 samples on the left and 500,000 on the right.}
\label{fig:prime_ran;vary_ss}
\end{figure}

To do this, we first generated a random sample of size $S$. Then, we replaced the last $x$\% of our sample, with a randomly chosen set of primes, sampled uniformly at random
from the first 50 million primes. We note that, in these tests, the originally randomly generated sample included both even and odd numbers.
From here, we calculated the Z-score for the our combined list. For simplicity, we started with a fixed sample size of 300,000 elements.\footnote{We also ran similar tests with sample sizes of 100,000 and 500,000. These graphs are in Figure \ref{fig:prime_ran;vary_ss}.} We varied our percentage of primes from 10\%
to 90\%. For each fixed percentage, we ran 1,000 trials. Our results are summarized in Figure \ref{fig:prime_ran;ss:300k}.

The relationship between Z-scores and the percentage of primes appears to be linear. We used a linear regression package to find the line of best fit. 
We found that the line $y = mx + b$, where $m = 0.221$ and $b=-0.071$ was the best fit line. This provided us with a correlation coefficient ($r^2$ value) of 0.987, which implies that our data is
highly likely to be linear. The p-value for the test where the null hypothesis is that the data is not linearly related (has a slope of zero) was $8.86 \times 10^{-7}$.

We can see that, even in the smaller sample sizes, there is still a linear relationship. The correlation coefficient for 100,000 samples and 500,000 samples are 0.9834 and 0.9887, respectively.
These produce p-values of $1.917\times 10^{-6}$ and $4.877 \times 10^{-7}$. Both of these results show a strong linear correlation between Z-score and the percentage of primes.

This shows that we can use the metric of the parity of the {\sc sum-of-digits} to approximate the percentage of prime numbers in a set of randomly generated numbers. That is, we can use our results,
based on the linear regression, to estimate how many more prime numbers there are in a set, compared to how many there would normally be in a set of random numbers.
This implies that our metric is not only useful for distinguishing between sets of random numbers; rather, it can also be used to estimate how ``prime'' a set of numbers is.

\section{Additional Tests}\label{sec:additional}

We went through several bouts of testing to confirm that the disparity in parities for primes was real. This included testing for bugs in our code and verifying our results with other people who independently reporduced the same results. We also tested changing the base, in case the results were unique to base 10. We tested our results on some other even bases as well. We note that, in odd bases, all prime numbers (other than 2) must have an odd sum of digits. If not, then, it has an even number of odd digits, which means that it corresponds to an even number in base 10, which implies that it can't be a prime, other than 2. Our results persisted in other bases, as well.

We also note that some of our experiments were run comparing prime numbers to a baseline of only odd random numbers, while others were comparing to a baseline of all random numbers. In base 10, having an odd last digit simply reverses the parity and shouldn't affect our results. We ran some sanity checks to test that these baselines behaved the same way, which they did.

We also ran some tests to see if the bias in the primes holds over modular arithmetic. In particular, we chose 100 prime numbers uniformly at random from the first 1,000 prime numbers.
Then, we chose $10^6$ random numbers. For each of these possible pairs, we computed $p \mod r$, where $p$ was the prime and $r$ was the random number. This gave us $10^8$ possible values,
which we expected to follow a binomial distribution with $10^8$ draws and $p=0.5$. However, in this test, we found that the prime numbers modded with random numbers,
showed a z-score of 14.18 deviations from the expected value. This shows that the bias in the prime numbers holds over modular arithmetic, which has possible repercussions in fields like cryptography.

\begin{figure}[!htbp]
\centering
 \includegraphics[width=10cm, height=6cm,keepaspectratio]{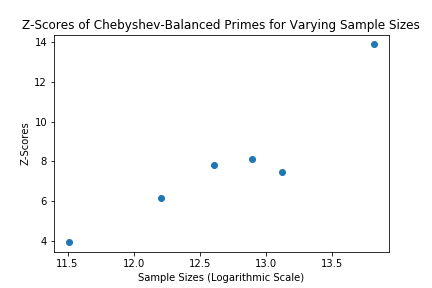}

\caption{Z-scores for Chebyshev-balanced primes, with varying sample sizes. Note that
the sample size is plotted on a logarithmic scale.}
\label{fig:chebyshev}
\end{figure}

The last test we ran was investigating whether the bias in the prime numbers was due to other inherent biases in the prime numbers. In particular, we tested Chebyshev's bias,
where primes are more likely to be in the form $(3 \mod 4)$, rather than $(1 \mod 4)$ \cite{chebyshev_bias}. To account for this, we created sets of prime numbers, which were split evenly between primes of the form $(3\mod 4)$ and $(1\mod 4)$, and tested if the sum of digits bias still held on similar sample sizes. Our results are summarized in Figure \ref{fig:chebyshev}.

We ran 100 trials for each test. In a given trial, we created a sample of size $s$ by sampling $s/2$ primes of the form $3 \mod 4$ and $s/2$ primes of the form $1 \mod 4$ from our list of the first 50 million primes.
Then, we computed the number of primes in our sample with an even sum of digits. We averaged our results over the 100 trials, and normalized to compute Z-scores, as before.
We can see that there is still a strong deviation from the expected value, even when the primes are evenly balanced, accounting to the Chebyshev bias.

Finally, we'd like to remark that the authors have spent some time investigating why the {\sc sum-of-digits} would provide such a disparity between prime numbers and random numbers. Our results indicate that the disparity is strong and appears to be a property of prime numbers. However, unlike metrics like prime number density, which can be justified theoretically, we have not been able to find a reason for this bias to occur for primes.

\section{Conclusion}\label{sec:conclusion}

We have described a metric, based on the parity of the {\sc sum-of-digits}, which can be used to distinguish between sets of randomly chosen prime numbers and sets of random numbers.

Using this metric, we have designed testing algorithms which succeed in distinguishing between sets of prime numbers and random numbers with high probability. We can also differentiate between a set of products-of-two-primes and a set of random numbers with probability greater than 50\%. We also conducted a number of experiments to show that the bias in the prime products is not solely due to some underlying bias of the sum-of-digits metric in the source set of prime numbers, but rather a intrinsic property of the prime numbers (and products) themselves.

Moreover, we tested mixed sets of prime and random numbers and showed a linear relationship between the percentage of primes added into a list of random numbers and the Z-score of the parity of the sum of digits. 

Finally, we ran several other tests including testing the primes modulo random numbers and balancing our samples of primes to account for Chebyshev's bias, to investigate causes of the {\sc sum-of-digits} bias. We found that the bias holds, even modulo random numbers but have not been able to find the root cause.

Future work will include many, many more tests along the lines of what we have conducted, with larger sample sizes, other metrics than Z-scores, comparisons between bases, etc. Of course, the most important thing would be to find a theoretical explanation for the bias; we have not had the faintest luck in finding one.

\bibliographystyle{unsrt}
\bibliography{primes}

\end{document}